\documentclass[aip,cha,amsfonts,amsmath,amssymb
,reprint,
]{revtex4-1}

\usepackage{amsthm}
\usepackage[dvips]{graphicx} 
\usepackage{subfig}
\usepackage{psfrag} 

\newtheorem{theorem}{Theorem}[section]
\newtheorem{proposition}[theorem]{Proposition}
\newtheorem{conjecture}[theorem]{Conjecture}
\newtheorem{corollary}[theorem]{Corollary}

\newtheorem{lemma}[theorem]{Lemma}

\newcommand{\Rr}{{\mathbb{R}}}

\newsavebox{\tempa}
\newsavebox{\tempb}

\begin{document}
\title
{Chaos in the square billiard with a modified reflection law}
\date{\today}

\author{Gianluigi Del Magno}
\email{delmagno@iseg.utl.pt}
\affiliation{CEMAPRE, ISEG, Universidade T\'ecnica de Lisboa, Rua do Quelhas 6, 1200-781 Lisboa, Portugal}
\author{Jo\~ao Lopes Dias}%
\email{jldias@iseg.utl.pt}
\affiliation{Departamento de Matem\'atica and CEMAPRE, ISEG, Universidade T\'ecnica de Lisboa, Rua do Quelhas 6, 1200-781 Lisboa, Portugal}
\author{Pedro Duarte}
\email{pduarte@ptmat.fc.ul.pt}
\affiliation{Departamento de Matem\'atica and CMAF, Faculdade de Ci\^encias, Universidade de Lisboa, Campo Grande, Edif\'icio C6, Piso 2, 1749-016 Lisboa, Portugal}
\author{Jos\'e Pedro Gaiv\~ao}
\email{jpgaivao@iseg.utl.pt}
\author{Diogo Pinheiro}
\email{dpinheiro@iseg.utl.pt}
\affiliation{CEMAPRE, ISEG, Universidade T\'ecnica de Lisboa, Rua do Quelhas 6, 1200-781 Lisboa, Portugal}

\begin{abstract}
The purpose of this paper is to study the dynamics of a square billiard with a non-standard reflection law such that the angle of reflection of the  particle is a linear contraction of the angle of incidence. We present numerical and analytical arguments that the nonwandering set of this billiard decomposes into three invariant sets, a parabolic attractor, a chaotic attractor and a set consisting of several horseshoes. This scenario implies the positivity of the topological entropy of the billiard, a property that is in sharp contrast with the integrability of the square billiard with the standard reflection law.
\end{abstract}

\pacs{05.45.Pq Numerical simulations of chaotic systems}
                            
\keywords{Billiard; Hyperbolicity; Horseshoe; Chaotic attractor}

\maketitle

\begin{quotation}
A billiard is a mechanical system consisting of a point-particle moving freely inside a planar region and being reflected off the perimeter of the region according to some reflection law. The specular reflection law is the familiar rule that prescribes the equality of the angles of incidence and reflection. Billiards with this reflection law are conservative systems, and as such are models for physical systems with elastic collisions. For this reason and their intrinsic mathematical interest, conservative billiards have been extensively studied. Much less studied are dissipative billiards, which originate from reflection laws requiring that the angle of reflection is a contraction of the angle of incidence. These billiards do not preserve the Liouville measure, and therefore can model physical systems with non-elastic collisions.
In this paper, we investigate numerically and analytically a dissipative billiard in a square. We find that its dynamics differs strikingly from the one of its conservative counterpart, which is well known to be integrable. Indeed, our results show that a dissipative billiard in a square has a rich dynamics with horseshoes and attractors of parabolic and hyperbolic type coexisting simultaneously.

\end{quotation}

\section{Introduction}

Billiards are among the most studied dynamical systems for two main reasons. Firstly, they serve as models for important physical systems (see e.g.~the book\cite{Szasz} and references therein), and secondly, despite their simplicity, they can display a rich variety of dynamics ranging from integrability to complete chaoticity. Most of the existing literature on billiards is devoted to billiards with the standard reflection law: the angle of reflection of the particle equals the angle of incidence (cf.~\cite{cm06,ta}). These billiards are conservative systems, i.e. they admit an invariant measure that is absolutely continuous with respect to the phase space volume.

In this paper, we are concerned with billiards with a non-standard reflection law according to which the angle of reflection equals the angle of incidence times a constant factor $ 0 < \lambda < 1 $. Since we have observed numerically that such a law has the effect of contracting the phase space volume, billiards with this law will be called ``dissipative'' in this paper.

Recently, Markarian, Pujals and Sambarino\cite{markarian10} proved that dissipative planar billiards (called ``pinball billiards'' in their paper)  
have two invariant directions such that the growth rate along one direction dominates uniformly the growth rate along the other direction. This property is called \emph{dominated splitting}, and is weaker than hyperbolicity, which requires one growth rate to be greater than one, and the other one to be smaller than one. The result of Markarian, Pujals and Sambarino applies to billiards in regions of different shapes. In particular, it applies to billiards in polygons. This is an interesting fact because the dominated splitting property enjoyed by the dissipative polygonal billiards contrasts with the parabolic dynamics observed in the conservative case \cite{markarian10,ta}.

Here, we take a further step towards the study of dissipative polygonal billiards analyzing the dissipative square billiard. Taking into account the symmetries of the square, we perform our analysis on a reduced phase space. We provide theoretical arguments and numerical evidence that the nonwandering set of our system decomposes into three (possibly empty) invariant sets: a parabolic attractor, a hyperbolic attractor and a horseshoe. This dynamics is clearly richer than the one of the conservative square billiard, which is a fully integrable system. 
In this paper, we also conduct a rather detailed numerical study of the changes in the properties of the nonwandering set as the parameter $ \lambda $ varies. 

We should mention that results somewhat similar to ours were obtained for non-polygonal billiards \cite{Altmann08,Arroyo_Markarian_Sanders} and the dissipative equilateral triangle billiard \cite{ams}.

The paper is organized as follows. In Section~\ref{se:square}, we give a detailed description of the map for the dissipative square billiard. Some results concerning the invariant sets of this map are presented in Section~\ref{se:hyperbolicity}. To study our map, it is convenient to quotient it by the symmetries of the square. This procedure is described in Section~\ref{se:reduced}, and produces the so-called reduced billiard map. Section~\ref{se:a-h} is devoted to the study of two families of periodic points of the reduced billiard map. In particular, we show the stable and unstable manifolds of a fixed point of the reduced billiard map (corresponding to a special periodic orbit of the billiard map) have transversal homoclinic intersections, and use this fact to conclude that the dissipative square billiard has positive topological entropy. Finally, Section \ref{se:bifurcation} contains the bifurcation analysis of the nonwandering set of the reduced billiard map. 

\section{The square billiard}
\label{se:square}

\begin{figure}[tb]
\includegraphics[width=0.45\textwidth]{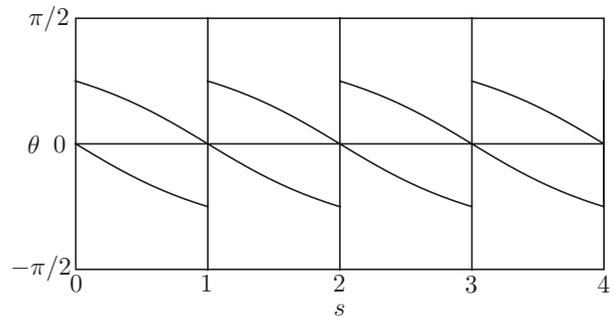}
\caption{Phase space $ \mathcal{M} \setminus \mathcal{S}^+$.}
\label{phase_space_fig}
\end{figure}

Consider the square $D=[0,1]\times[0,1]\subset\Rr^2$. For our purposes, $D$ is called the square billiard table. 
To study the dynamics of the billiard inside this table, it is sufficient to know the angle of incidence at the impact points and the reflection law. For the usual reflection law (the angle of reflection is equal to the angle of incidence) we find the next impact point $s'$ and angle of reflection $\theta'$ by the billiard map $(s',\theta')=\mathcal{B}(s,\theta)$ acting on the previous impact $(s,\theta)$. This map admits an explicit analytic description. Its domain coincides with the rectangle  
$$
\mathcal{M}=[0,4]\times\left(-\frac\pi2,\frac\pi2\right)
$$
from which the set
$$
\mathcal{S}^+=\left\{\, (s,\theta)\in \mathcal{M} \colon \{s\} = 0 \text{ or } \{s\}+\tan\theta \in \{0,1\} \, \right\}
$$
is removed. The symbols $[s]$ and $\{s\}=s-[s]$ stand for the integer part and the fractional part of $s$, respectively. An element of $ \mathcal{S}^{+} $ corresponds to an orbit leaving or reaching a corner of $ D $ (see Fig.~\ref{phase_space_fig}). 
By reversing the role of time in this description of $ \mathcal{S}^{+} $, one obtains the set
$$
\mathcal{S}^-=\left\{\, (s,\theta)\in \mathcal{M} \colon \{s\} = 0 \text{ or } (s,-\lambda^{-1} \theta) \in \mathcal{S}^{+} \, \right\}.
$$
Both sets $ \mathcal{S}^{+} $ and $ \mathcal{S}^{-} $ consist of finitely many analytic curves. Next, let
\begin{align*}
\mathcal{M}_{1} & = \left\{\,(s,\theta) \in \mathcal{M} \colon \{s\} > 0 \text{ and } \{s\}+\tan\theta > 1 \,\right\}, \\
\mathcal{M}_{2} & = \left\{\,(s,\theta) \in \mathcal{M} \colon \{s\} > 0 \text{ and } 0 < \{s\}+\tan\theta < 1 \,\right\}, \\
\mathcal{M}_{3} & = \left\{\,(s,\theta) \in \mathcal{M} \colon \{s\} > 0 \text{ and } \{s\}+\tan\theta < 0 \,\right\}.
\end{align*}
The billiard map $ \mathcal{B} \colon \mathcal{M} \setminus \mathcal{S}^{+} \to  \mathcal{M} \setminus \mathcal{S}^{-} $ is defined by
$$
\mathcal{B}(s,\theta)=
\begin{cases}
\left([s]+1+\dfrac{1-\{s\}}{\tan \theta} \negthickspace \negthickspace \pmod 4,\dfrac{\pi}{2} - \theta\right) & \negthickspace \negthickspace \text{on } \mathcal{M}_{1}, \\ 
\left([s]-1-\{s\}-\tan\theta \negthickspace \negthickspace \pmod 4, -\theta\right) & \negthickspace \negthickspace \text{on } \mathcal{M}_{2}, \\
\left([s]+\dfrac{\{s\}}{\tan\theta} \negthickspace \negthickspace \pmod 4,-\dfrac\pi2-\theta \right) 
& \negthickspace \negthickspace \text{on } \mathcal{M}_{3}.
\end{cases}
$$
This map is clearly an analytic diffeomorphism in its domain.
The inverse of $\mathcal{B}$ is easily obtained by noticing that the billiard map is time-reversible.
That is, given the map $\mathcal{T}(s,\theta)=(s,-\theta)$, we have
$$
\mathcal{B}^{-1}=\mathcal{T}\circ \mathcal{B}\circ\mathcal{T}^{-1}.
$$

To modify the reflection law, we compose $ \mathcal{B} $ with another map $\mathcal{R} \colon \mathcal{M} \to \mathcal{M}$. The resulting map $ \Phi = \mathcal{R} \circ \mathcal{B} $ is called a \emph{billiard map with a modified reflection law}.   

\begin{figure}[t]
\includegraphics[width=0.2\textwidth]{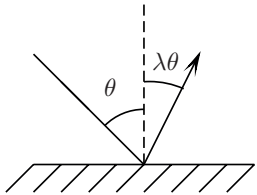}
\caption{Dissipative reflection law.}
\label{reflection_law}
\end{figure}

Several reflections laws have been considered~\cite{Altmann08,markarian10}.
In this paper, we consider the following ``dissipative'' law. Given $0<\lambda<1$, we set
$$
\mathcal{R}_\lambda(s,\theta)=\left(s,\lambda\theta\right).
$$
According to this law, the direction of motion of the particle after a reflection gets closer to the normal of the perimeter of the square (see Fig.~\ref{reflection_law}). To emphasize the dependence of the billiard map on the parameter $ \lambda $, we write 
$$
\Phi_\lambda=\mathcal{R}_\lambda\circ \mathcal{B}.
$$

As a side remark, one can also define the map $ \Phi_{\lambda} $ for $ \lambda > 1 $. In this case, the map $ \mathcal{R}_{\lambda} $ expands uniformly the angle $ \theta $, and $ \Phi_{\lambda} $ becomes a map with holes in the phase space. It is interesting to observe that the maps $ \Phi_{\lambda^{-1}} $ and $ \Phi^{-1}_{\lambda} $ are conjugated for $ 0<\lambda<1$. Indeed, it is not difficult to check that 
\[
\Phi_{\lambda^{-1}} = (\mathcal{R}_{\lambda}\circ \mathcal{T})^{-1} \circ \Phi^{-1}_{\lambda} \circ (\mathcal{R}_{\lambda}\circ \mathcal{T}),
\]
by using the fact that $\mathcal{T}$ and $\mathcal{R}_{\lambda}$ commute and that $\mathcal{R}_{\lambda}^{-1}=\mathcal{R}_{\lambda^{-1}}$.
Therefore, all the results presented in this paper hold for $ \lambda > 1 $ as well, provided that we replace the word ``attractor'' with the word ``repeller'', and switch the words ``stable'' and ``unstable''. 

\section{Hyperbolicity} 
\label{se:hyperbolicity}

Let $ (s_{0},\theta_{0}) $ be an element of $ \mathcal{M} \setminus \mathcal{S}^{+} $. Set $ (s_{1},\theta_{1}) = \Phi_{\lambda}(s_{0},\theta_{0}) $, and denote by $ t(s_{0},\theta_{0}) $ the length of the segment connecting $ s_{0} $ and $ s_{1} $. 
Using elementary trigonometry, one can show in a straightforward manner that the derivative of $ \Phi_{\lambda} $ takes the following form:
\begin{equation*}
	D \Phi_{\lambda}(s_{0},\theta_{0}) = -
	\begin{pmatrix}
		\dfrac{\cos \theta_{0}}{\cos \lambda^{-1} \theta_{1}} & \dfrac{t(s_{0},\theta_{0})}{\cos \lambda^{-1} \theta_{1}} \\
		0 & \lambda \\
	\end{pmatrix}.
\end{equation*}
In fact, the previous formula holds for every polygon, and not just for the square (see Formula 2.26 in ref.~\cite{cm06}).

Now, suppose that $ \{(s_{i},\theta_{i})\}^{n}_{i=0} $ are $ n+1 $ consecutive iterates of $ \Phi_{\lambda} $.
Then, we see that
\begin{equation*}
	D\Phi_{\lambda}^{n} (s_{0},\theta_{0}) = (-1)^{n}
	 \begin{pmatrix}
			\alpha_{n}(s_{0},\theta_{0})  & \zeta_{n}(s_{0},\theta_{0}) \\
			0 & \beta_{n}(s_{0},\theta_{0}) \\
	\end{pmatrix},
\end{equation*}
where
\[
\alpha_{n}(s_{0},\theta_{0}) = \frac{\cos \theta_{0}}{\cos \lambda^{-1} \theta}_{n} \prod^{n-1}_{i=1} \frac{\cos \theta_{i}}{\cos \lambda^{-1} \theta}_{i}, \quad \beta_{n}(s_{0},\theta_{0}) = \lambda^{n},
\]
and
\[ 
\zeta_{n}(s_{0},\theta_{0}) = \frac{1}{\cos \lambda^{-1} \theta_{n}} \sum^{n-1}_{i=0} \lambda^{i} t(s_{i},\theta_{i}) \prod^{n-1}_{k=i+1} \frac{\cos \theta_{k}}{\cos \lambda^{-1} \theta}_{k}.
\] 

We now prove a simple lemma concerning the stability of the periodic points of $ \Phi_{\lambda} $. It is not difficult to see that this result remains valid for every polygon and for other reflection laws (e.g. $ \mathcal{R}_{\lambda}(s,\theta) = (s,\theta - c \sin 2 \theta) $ with $ 0 < c < 1/2 $ as in~\cite{markarian10}). 
	
\begin{lemma}
\label{le:periodic-points}
For every $ \lambda \in (0,1) $, the periodic points of $ \Phi_{\lambda} $ of period 2 and period greater than 2 are parabolic and hyperbolic, respectively.
\end{lemma}
	
\begin{proof}	
	Suppose that $ (s_{0},\theta_{0}) $ is a periodic point of $ \Phi_{\lambda} $ with period $ n $. Since $ (s_{n},\theta_{n}) =  (s_{0},\theta_{0}) $, it turns out that 
\[ 
\alpha_{n}(s_{0},\theta_{0}) = \prod^{n-1}_{i=0} \frac{\cos \theta_{i}}{\cos \lambda^{-1} \theta_{i}}.
\] 

Now, note that each term $ \cos \theta_{i} / \cos \lambda^{-1} \theta_{i} $ in the expression of $ \alpha_{n}(s_{0},\theta_{0}) $ is equal or greater than 1 with equality if and only if $ \theta_{i} = 0 $. Also, note that $ D\Phi_{\lambda}^{n} (s_{0},\theta_{0}) $ is a triangular matrix, and so the moduli of its eigenvalues are $ \alpha_{n}(s_{0},\theta_{0}) $ and $ \lambda^{n} < 1  $. Therefore to determine the stability of $ (s_{0},\theta_{0}) $ is enough to find out whether or not $ \alpha_{n}(s_{0},\theta_{0}) $ is greater than 1.

If $ n = 2 $, it is easy to see that the trajectory of $ (s_{0},\theta_{0}) $ must always hit the boundary of $ D $ perpendicularly. In other words, we have $ \theta_{0} = \theta_{1} = \theta_{2} = 0 $, and so $ \alpha_{2}(s_{0},\theta_{0}) = 1 $. Periodic points of period 2 are therefore parabolic. Clearly, a necessary condition for a polygon to admit periodic points of period 2 is that the polygon must have at least 2 parallel sides (not a sufficient condition though). 

Now, suppose that $ n > 2 $. In this case, we claim that $ (s_{0},\theta_{0}) $ is hyperbolic. Indeed, when $ n > 2 $, the billiard trajectory of $ (s_{0},\theta_{0}) $ must have at least two non-perpendicular collisions with the boundary of $ D $, and since $ \cos \theta_{i}/ \cos \lambda^{-1} \theta_{i} > 1 $ for such collisions, we can immediately conclude that $ \alpha_{n}(s_{0},\theta_{0}) > 1 $. 
\end{proof}

A more elaborated analysis along the lines of the proof of Lemma \ref{le:periodic-points} yields some general conclusions on the chaotic behavior of general dissipative polygonal billiards. Two of such conclusions are stated in Propositions \ref{pr:delta-positive} and \ref{pr:hyperbolicity} below. To state these proposition, we need first to introduce the notion of uniformly and non-uniformly hyperbolic systems. Unless specified otherwise, $ \Phi_{\lambda} $ denotes the map of a dissipative billiard in a general polygon $ D $ throughout the rest of this section.

A set $ \Sigma \subset \mathcal{M} $ is called \emph{invariant} if $ \Phi^{-1}_{\lambda}(\Sigma) = \Sigma $. An invariant set $ \Sigma $ is called \emph{hyperbolic} if there exist a norm $ \|\cdot\| $ on $ \mathcal{M} $, a non-trivial invariant measurable splitting $ T_{\Sigma} \mathcal{M} = E^{s} \oplus E^{u} $ and two measurable functions $ 0 < \mu < 1 $ and $ K > 0 $ on $ \Sigma $ such that for every $ (s,\theta) \in \Sigma $ and every $ n \ge 1 $, we have
\begin{align*}
\|D\Phi^{n}_{\lambda}|_{E^{s}(s,\theta)}\| & \le K(s,\theta) \mu(s,\theta)^{n}, \\
\|D\Phi^{-n}_{\lambda}|_{E^{u}(\Phi^{n}_{\lambda} (s,\theta))}\| & \le K(s,\theta) \mu(s,\theta)^{n}.
\end{align*}
If the functions $ \mu $ and $ K $ can be replaced by constants, then $ \Sigma $ is called \emph{uniformly hyperbolic}, otherwise it is called \emph{non-uniformly hyperbolic}. 

We can now state our propositions. The first one concerns billiards in polygons without parallel sides. For such polygons, the map $ \Phi_{\lambda} $ does not have periodic points of period 2. 

\begin{proposition}
	\label{pr:delta-positive}
	Let $ D $ be a polygon without parallel sides, and suppose that $\Sigma$ is an invariant set of $ \Phi_{\lambda} $. Then $ \Sigma $ is uniformly hyperbolic for every $ \lambda \in (0,1) $.
\end{proposition}

The second proposition concerns billiards in rectangles. In this case, $ \Phi_{\lambda} $ has periodic points of period 2. Denote by $ \mathcal{P} $ the set of all these points. It is not difficult to check that $ \mathcal{P} $ is a parabolic attractor for every $ \lambda \in (0,1) $. 

\begin{proposition}
	\label{pr:hyperbolicity}
	Let $ D $ be a rectangle, and suppose that $\Sigma$ is an invariant set of $ \Phi_{\lambda} $ not intersecting the basin of attraction of $ \mathcal{P} $. Then there exists $ \lambda_{*} \in (0,1) $ such that $ \Sigma $ is hyperbolic for every $ \lambda \in (0,\lambda_{*}) $, and is uniformly hyperbolic for every $ \lambda \in (\lambda_{*},1) $.
\end{proposition}

For reasons of space, the proofs of these propositions are omitted and will appear elsewhere~\cite{MDDGP}.

\section{The reduced billiard map}
\label{se:reduced}

The analysis of the billiard dynamics can be simplified if we reduce the phase space. First, we identify all sides of the square by taking the quotient with the translations by integers of the $s$-component. Then, due to the symmetry along the vertical axis at the midpoint of the square, we can also identify each point $(s,\theta)$ with $(1-s,-\theta) $. To formulate the reducing procedure more precisely, we define an equivalence relation $\sim$ on $\mathcal{M}$ by  
$(s_1,\theta_1)\sim(s_2,\theta_2)$ if and only if  
$\pi(s_1,\theta_1)=\pi(s_2,\theta_2)$, where $\pi\colon \mathcal{M}\to \mathcal{M}$  
is the function defined by
\begin{equation*}
\pi(s,\theta) =
\begin{cases}
(\{s\}, \theta) & \text{if } \theta \in \left[0,\frac\pi2\right), \\
(1-\{s\}, -\theta) & \text{if } \theta \in \left(-\frac\pi2,0\right).
\end{cases}
\end{equation*}
Let $M$ denote the image of $\pi$. Clearly, we have 
$$
M=(0,1)\times\left[0,\frac\pi2\right).
$$
Note that it is possible to identify the set $M$ with the quotient space $\mathcal{M} /\sim$. We call $M$ the \emph{reduced phase space}. The induced billiard map on $M$ is the \emph{reduced map}, which we will denote by $\phi_\lambda$.

It is clear from the definition of $\pi$ that $ \pi^{-1}(s,\theta)$ consists of 8 elements for every $(s,\theta)\in M$, and so $(\mathcal{M},\pi)$ is an 8-fold covering of $M$. It is then easy to see that the reduced billiard map $\phi_\lambda$ is a factor of the original billiard map $\Phi_\lambda$ by noting that the quotient map $\pi$ is indeed a semiconjugacy between $\phi_\lambda$ and $\Phi_\lambda$, i.e. we have that $\pi\circ\Phi_\lambda=\phi_\lambda \circ \pi$. 

In what concerns the relation between the dynamical systems defined by $\Phi_\lambda$ and $\phi_\lambda$, there are several key points that are worth remarking. First, we note that periodic points of $\phi_\lambda$ lift to periodic points of $\Phi_\lambda$. To be more precise, an orbit of period $n$ under $\phi_\lambda$ is lifted to either eight orbits of period $n$, or four orbits of period $2n$, or two orbits of period $4n$, or one orbit of period $8n$ for $\Phi_\lambda$. Analogous statements can be produced for the lifts of transitive sets and the existence of invariant measures. Namely, transitive sets for $\phi_\lambda$ are lifted to a finite number of transitive sets for $\Phi_\lambda$, and any invariant measure under the dynamics of $\phi_\lambda$ corresponds to a finite number of invariant measures under $\Phi_\lambda$. Finally, we remark that the reduced map $\phi_\lambda$ has positive topological entropy if and only if this is the case for the billiard map $\Phi_\lambda$.

By studying the trajectories of the billiard map we have basically two cases: either the billiard orbit hits a neighboring side of the square or the opposite side.
Separating these cases there is a corner which is reachable only if the initial position $(s,\theta)\in M$ is in the singular curve
$$
S^+=\left\{
(s,\theta)\in M\colon
s+\tan\theta=1
\right\}.
$$
This curve separates the reduced phase space in two connected open sets: $ M_1$ below $S^+$ and $ M_2$ above $ S^+ $. 

Let $f_1 \colon M_{1} \to M $ and $f_2 \colon M_{2} \to M $ be the transformations defined by
\begin{align*}
f_1(s,\theta) & = \left(s+\tan\theta,\lambda\theta\right) & \text{for } (s,\theta) \in M_{1}, \\
f_2(s,\theta) & = \left((1-s) \cot \theta,\lambda\left(\frac\pi2-\theta\right)\right) & \text{for } (s,\theta) \in M_{2}.
\end{align*}
The reduced billiard map is then given by
$$
\phi_\lambda =
\begin{cases}
	f_1 & \text{on } M_1, \\
	f_2 & \text{on } M_2. 
\end{cases}
$$
Its domain and range are $M\setminus S^+$ and $ M \setminus S^{-} $, respectively, where 
$$
S^-=\left\{
(s,\theta)\in M\colon
s-\tan(\lambda^{-1}\theta)=0
\right\}.
$$
Like the billiard map $ \Phi_{\lambda} $, the reduced billiard map $ \phi_{\lambda} $ is an analytic diffeomorphism.
Notice that $\phi_\lambda$ maps horizontal lines into horizontal lines, a consequence of the fact that its second component is independent of $s$.

Finally, we observe that the subsets of $ M $ where the maps $ \phi^{n}_{\lambda} $ and $ \phi^{-n}_{\lambda} $ are defined for every $ n \ge 0 $ are, respectively,
$$
M^{+}=M \setminus \bigcup_{n \ge 0} \phi^{-n}_{\lambda}(S^{+}) \quad \text{and} \quad M^{-}=M \setminus \bigcup_{n \ge 0} \phi^{n}_{\lambda}(S^{-}).
$$


\begin{figure}[tb]
\includegraphics[width=0.28\textwidth]{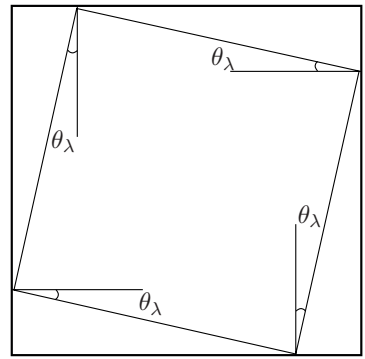}
\caption{Periodic orbit.}
\label{fig:periodic}
\end{figure}

\section{Attractors and horseshoes}
\label{se:a-h}

We start this section by formulating several definitions. 

The {\em stable set} of an element $ q \in M $ is defined by 
$$
W^s(q)=\left\{\, x\in  M^{+}\colon \lim_{n\to+\infty} \left\| \phi_\lambda^n(x)-\phi_\lambda^n(q) \right\| = 0 \,\right\},
$$
where $ \|\cdot\| $ is the Euclidean norm on $ M $. In the case of an invariant set $\Lambda=\phi_\lambda(\Lambda)$, we define its stable set to be 
$$
W^s(\Lambda)=\bigcup_{q\in\Lambda}W^s(q).
$$
The {\em unstable sets} $W^u(q)$ and $W^u(\Lambda)$ are defined analogously by replacing $ \phi_{\lambda} $ with $ \phi^{-1}_{\lambda} $ and $ M^{+} $ with $ M^{-} $. When $W^{u(s)}(\Lambda)$ turns out to be a manifold, we will call it an {\em unstable}({\em stable}) manifold.

Suppose that $\Lambda$ is an invariant subset of $M$. Then $\Lambda$ is called an {\em attractor} if $\Lambda=W^u(\Lambda)$ and $W^s(\Lambda)$ is open in $M^{+}$, and is called a {\em horseshoe} if neither $W^s(\Lambda)$ is an open set in $M^{+}$ nor is $W^u(\Lambda)$ an open set in $M^{-}$. Note that a saddle periodic orbit is a horseshoe according to this definition. A finite union of hyperbolic invariant sets $A_1,\ldots, A_m$ is called a {\em hyperbolic chain} if
$$
W^u(A_i)\cap W^s(A_{i+1})\neq \emptyset \qquad \text{for } i=1,\ldots,m-1. 
$$

A point $x\in M_+$ is said to be {\em nonwandering} if for every open neighborhood $U$ containing $x$, there exists $n\geq 1$ such that $U\cap \phi_\lambda^ {n}(U)\neq \emptyset$. We denote by $\Omega_{\lambda} \subset M^ +\cap M^-$ the set of all nonwandering points of $\phi_{\lambda}$. We say that two hyperbolic periodic points  $x,y\in \Omega_{\lambda}$ are {\em homoclinically related} if $W^u(x)$ and $W^s(y)$ intersect transversally, and $W^u(y)$ and $W^s(x)$ intersect transversally. The closure in $M^+\cap M^-$ of the set of periodic points homoclinically related to a hyperbolic periodic point $x \in \Omega_{\lambda} $ is called the {\em homoclinic class of $x$}. Every homoclinic class is a transitive invariant subset of $\Omega_{\lambda}$ (see Ref.~\cite[Ch.~IX, Prop.~5.2]{r}).

\subsection{Parabolic attractor}\label{ss:parabolic_attractor}

Let us define $P=\pi(\mathcal{P})$. It is easy to see that
$$
P=\left\{(s,\theta)\in M\colon \theta=0\right\},
$$
and each point of $P$ is a parabolic fixed point coming from period 2 orbits of the original billiard (orbits that bounce between parallel sides of the square). It is an attractor and $ W^{s}(P) $ includes the set of points $B$ that are below the forward invariant curve
$$
S_{\infty}=\left\{
(s,\theta)\in M\colon
s+\sum_{i=0}^{+\infty}\tan(\lambda^i\theta)=1
\right\}.
$$
The sequence $\phi_\lambda^n(S_{\infty})$ converges to the point $(1,0)$. The pre-image of $B$ is at the top of phase space. Moreover, its basin of attraction is
$$
W^s(P)=\bigcup_{n\geq0} \phi_\lambda^{-n}(B).
$$

By Proposition~\ref{pr:hyperbolicity}, the set $\Omega_{\lambda}\setminus P$ is hyperbolic. Hence, every periodic point in $\Omega_{\lambda} \setminus P$ has stable and unstable manifolds. Because of the cutting effect of the singular sets $S^{-}$ and $S^{+}$, these manifolds are countable unions of smooth curves. 

\begin{figure}[t]
\includegraphics[width=0.3\textwidth]{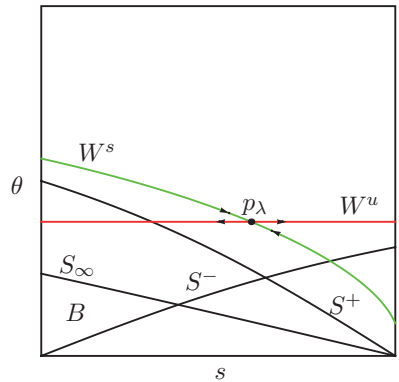}
\caption{Invariant manifolds of $p_\lambda$ and singular curves for the reduced billiard map ($\lambda=0.6218$)}
\label{reduced_phase_space_fig}
\end{figure}

\subsection{Fixed point and its invariant manifolds}

The map $ \Phi_{\lambda} $ has many periodic orbits. Two special periodic orbits of period 4 can be found by using the following simple argument. A simple computation shows that if an orbit hits two adjacent sides of the square with the same reflection angle $\theta_{\lambda}$, then 
$$
\theta_{\lambda}=\frac{\pi\lambda}{2(1+\lambda)}.
$$
If we further impose the condition that the orbit hits the two sides at $s_{1}$ and $s_{2}$ in such a way that $\{s_{1}\}=\{s_{2}\}=s_{\lambda}$, then we obtain
$$
s_{\lambda}=\frac1{1+\tan\theta_{\lambda}}.
$$
By symmetry, we conclude that $\Phi_\lambda^4(s_{\lambda},\theta_{\lambda})=(s_{\lambda},\theta_{\lambda})$. Using once again the symmetry of the square, we also have $\Phi^{4}_{\lambda}(1-s_{\lambda},-\theta_{\lambda})=(1-s_{\lambda},-\theta_{\lambda})$. One of these orbits is depicted in Fig.~\ref{fig:periodic}.

Due to the phase space reduction, the periodic orbits just described correspond to the fixed point 
\[ 
p_\lambda=(s_\lambda,\theta_\lambda) 
\] 
of $\phi_\lambda$. This is actually the only fixed point of $ \phi_{\lambda} $ in $M_2$ outside of $P$. By Lemma \ref{le:periodic-points}, $ p_{\lambda} $ is hyperbolic and thus it has local stable and unstable manifolds $W^{s,u}_{\rm loc}(p_\lambda)$ for every $\lambda\in(0,1)$. 
Since $\phi_\lambda$ maps horizontal lines into horizontal lines, and the set $ S^{-} $ does not intersect the horizontal line through $ p_{\lambda} $, we see that the local unstable manifold of $p_\lambda$ is given by
$$
W^u_{\rm loc}(p_\lambda)=\left\{\,(s,\theta)\in M\colon \theta=\theta_\lambda\,\right\}.
$$
In fact, the global unstable manifold consists of a collection of horizontal lines cut by the images of $S^{-}$. 

The geometry of the stable manifold is more complicated. By definition points on the stable manifold converge to the fixed point. Moreover, $W_{\rm loc}^s(p_\lambda)$ cannot cross $S^+$. Thus, $W_{\rm loc}^s(p_\lambda)$ is contained in $ M_2$. The graph transform associated with the corresponding branch of $\phi_\lambda$ is the transformation
$$
\Gamma(h)(\theta)= 1- h(g_\lambda(\theta))\,\tan\theta,
$$
where $g_\lambda\colon [0,\pi/2)\to [0,\pi/2)$ denotes the affine contraction
$$
g_\lambda(\theta)= \lambda\left(\frac{\pi}{2}-\theta\right).
$$
Iterating  $k$ times the zero function by $\Gamma$ we obtain
$$ 
\Gamma^k(0)(\theta) = \sum_{n=0}^{k-1} (-1)^n\,\prod_{i=0}^{n-1} \tan (g_\lambda^i(\theta) ) .
$$
Hence, the local stable manifold of $p_\lambda$ is the curve 
$$
W^s_{\rm loc}(p_\lambda)=\left\{\, (h_\lambda(\theta),\theta) \colon 0\leq\theta < \frac\pi2 \text{ and } 0 < h_\lambda(\theta) < 1 \, \right\},
$$ 
where
\begin{equation}\label{h:series}
h_\lambda(\theta) = \sum_{n=0}^{\infty} (-1)^n\,\prod_{i=0}^{n-1} \tan (g_\lambda^i(\theta)).
\end{equation}
This series converges uniformly and absolutely since $\tan(g_\lambda^n(\theta))$ converges to $ \tan \theta_\lambda $ as $n\rightarrow \infty$, and $0<\tan\theta_\lambda <1$. The same statement holds for the series of the derivatives of $h_\lambda$. Thus, $h_\lambda$ is smooth. 

The invariant manifolds of $p_\lambda$, the singular curves of the reduced billiard map and the upper boundary $S_\infty$ of $B$ are depicted in Fig.~\ref{reduced_phase_space_fig}. 

\begin{figure*}
 \centering
\sbox{\tempa}{\includegraphics[width=.41\textwidth]{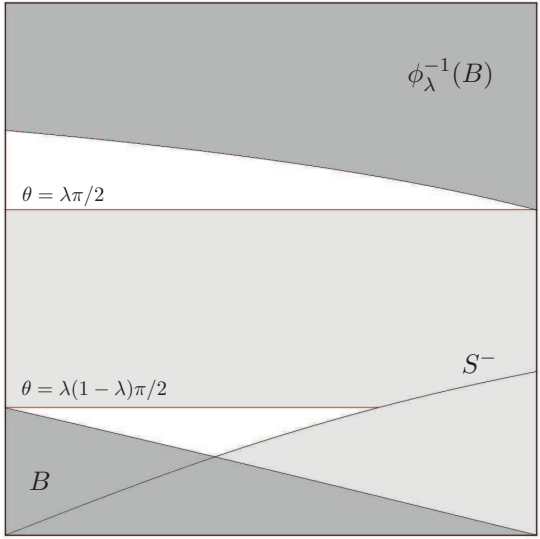}}%
\sbox{\tempb}{\includegraphics[width=.41\textwidth]{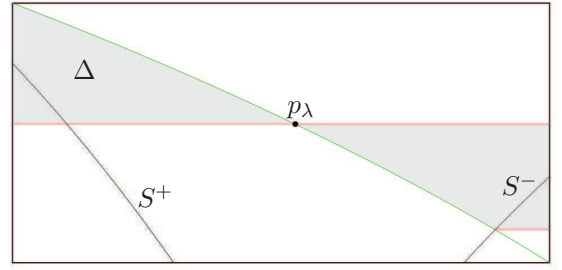}}%
\subfloat[$\lambda=\lambda_1$. Since the image of the map $\phi_\lambda$ is always below the line $\theta=\lambda\pi/2$, the basin of attraction of $P$ is only $B \cup \phi_\lambda^{-1}(B)$. Moreover, the region in light gray is forward-invariant.]{%
\begin{minipage}[b]{0.48\linewidth}\centering%
\usebox{\tempa}
\end{minipage}
}%
\hfill%
\subfloat[$\lambda=\lambda_2$. The shaded region $\Delta$ between the stable and unstable local manifolds of $p_\lambda$ is forward-invariant.]{%
\begin{minipage}[b]{0.48\linewidth}\centering%
\raisebox{1.5 cm}{\usebox{\tempb}}
\end{minipage}
}%
\caption{Trapping regions for $\phi_\lambda$.}
\label{trapping:region}
\end{figure*}

Let $ \lambda_{2}$ be the unique solution of
$$
h_{\lambda}(\lambda\theta_\lambda)=\tan(\theta_\lambda)\,\quad\text{for}\quad\lambda\in(0,1)\,.
$$
Geometrically, $ \lambda_{2} $ is the value of $ \lambda $ such that the singular set $ S^{-} $, the local manifold $W^s_{\mathrm{loc}}(p_\lambda)$ and the closure of the first iterate of the piece of the unstable manifold of $ p_{\lambda} $ contained in $ M_{1} $ have non-empty intersection (see Fig.~\ref{trapping:region}(b)). A numerical computation shows that 
$$
\lambda_2= 0.8736...
$$

Let $\Delta$ be the closed set bounded by $W^s_{\mathrm{loc}}(p_\lambda)$ and $W^u_{\mathrm{loc}}(p_\lambda)$ as in Fig.~\ref{trapping:region}(b).

\begin{proposition} \label{prop:delta}
For $\lambda>\lambda_2$, there is a compact $\phi_\lambda$-invariant set $\Delta_0\subset {\rm int}(\Delta)$ such that $\Omega_{\lambda}\cap \Delta \subset \{p_\lambda\}\cup \Delta_0$.
\end{proposition}
\begin{proof}
For $\lambda>\lambda_2$, we have $f_1(\Delta\cap M_1)\subset {\rm int}(\Delta)$, and we can a find a compact forward-invariant set $\Delta_0\subset {\rm int}(\Delta)$ under $ f_{2}$ such that $f_1(\Delta\cap M_1)\subset \Delta_0$. Hence $\Delta_0$ is also $\phi_\lambda$-invariant.
Since $p_\lambda$ is the only nonwandering point in $\Delta$ whose orbit does enter $M_1$, we have   $\Omega_{\lambda}\cap \Delta \subset \{p_\lambda\}\cup \Delta_0$.
\end{proof}

\begin{proposition}
	\label{pr:homoclinic}
The invariant manifolds of $ p_{\lambda} $ have transverse homoclinic points if and only if $ 0<\lambda<\lambda_{2} $.
\end{proposition}
\begin{proof}
To prove the existence of homoclinic points we iterate a piece of the local unstable manifold in $M_1$ and show that it intersects transversely the local stable manifold in $M_2$. Taking into account that $\phi_\lambda$ maps horizontal lines into horizontal lines, $\left.\phi_\lambda\right|_{M_1}=f_1$ and 
\[
f_{1}(W^{u}_{\textrm{loc}} \cap M_{1}) = (\tan \theta_{\lambda},1) \times \{\lambda \theta_{\lambda}\} 
\] 
the problem of finding homoclinic intersections reduces to proving the following chain of inequalities:
\begin{equation*}
\tan\theta_\lambda<h_\lambda(\lambda\theta_\lambda)<1.
\end{equation*}
We will see that these inequalities hold if and only if $0<\lambda<\lambda_2$. 
\begin{lemma}
The inequality $\tan\theta_{\lambda}<h_\lambda(\lambda\theta_{\lambda})$ holds if and only if $0<\lambda<\lambda_2$\,.
\end{lemma}
\begin{proof}
For $\lambda>0$ sufficiently small we have $\tan\theta_\lambda+\tan(\lambda\theta_\lambda)<1$ since $\theta_\lambda\rightarrow0$ as $\lambda\rightarrow0$. On the other hand, we know by definition of $h_\lambda$ that $h_\lambda(\lambda\theta_\lambda)+\tan(\lambda\theta_\lambda)>1$ for every $\lambda \in(0,1)$. Putting these two inequalities together we conclude that 
$$
\tan\theta_\lambda<1-\tan(\lambda\theta_\lambda)<h_\lambda(\lambda\theta_\lambda)
$$ 
for every $\lambda>0$ sufficiently small. Since $h_\lambda(\lambda\theta_\lambda)-\tan\theta_\lambda$ is strictly decreasing for $\lambda\in(0,1)$ and $\lambda_2$ is the unique solution of $h_\lambda(\lambda\theta_\lambda)=\tan\theta_\lambda$ we obtain the desired result. 
\end{proof}

\begin{lemma}\label{pr:h-leq-1}
The inequality $h_\lambda(\lambda\,\theta_\lambda)<1$ holds for every $\lambda\in (0,1)$.
\end{lemma}
\begin{proof}
At the fixed point we compute
\begin{equation}\label{eq:hprimefixpt}
h_\lambda'(\theta_\lambda)= -\frac{\sec^2\theta_\lambda}{(1-\lambda\,\tan \theta_\lambda)\,(1+\tan \theta_\lambda)} <0\;.
\end{equation}
Define now 
\begin{equation*}\begin{split}
m(\theta)&=h_\lambda(\theta_\lambda)+h_\lambda'(\theta_\lambda)\,(\theta-\theta_\lambda)\\
&=
 \frac{1-\lambda\,\tan \theta_\lambda -  \sec^2\theta_\lambda\,(\theta-\theta_\lambda)  }{ (1-\lambda\,\tan \theta_\lambda)\,(1+\tan \theta_\lambda) }\;
\end{split}
\end{equation*}
to be the function whose graph is the line tangent to the graph of $h_\lambda$ at the fixed point $p_\lambda$.
Since $h_\lambda$ is concave, $h_\lambda(\theta)\leq m(\theta)$ for every $\theta\in (0,\pi/2)$.
Thus, it is enough to check that $m(\lambda\,\theta_\lambda)< 1$ for $\lambda<1$, by using elementary estimates.
\end{proof}

Thus, $f_1(W^u_{\rm loc}(p_\lambda)\cap M_1)\cap W^s_{\rm loc}(p_\lambda)\not=\emptyset$ if and only if $0<\lambda<\lambda_2$. To conclude the proof of the proposition, note that for $\lambda>\lambda_2$, by Proposition~\ref{prop:delta} the region $\Delta$ is a trapping set. Since $W^u(p_\lambda)\subset\mathrm{int}(\Delta)\cup W^u_{\mathrm{loc}}(p_\lambda)$, $p_\lambda$ has no homoclinic intersections for $\lambda>\lambda_2$. 
\end{proof}

The following corollary is a direct consequence of Proposition~\ref{pr:homoclinic} and Ref.~\cite[Ch.~7,Th.~4.5]{r}.
\begin{corollary}
The map $\phi_\lambda$ has positive topological entropy for every $0<\lambda<\lambda_2$.
\end{corollary}

\subsection{Two Families of Periodic Orbits }\label{sect Two Families of Periodic Orbits}

\begin{figure*}
\subfloat[$q_1$]{\includegraphics[width=.2\textwidth]{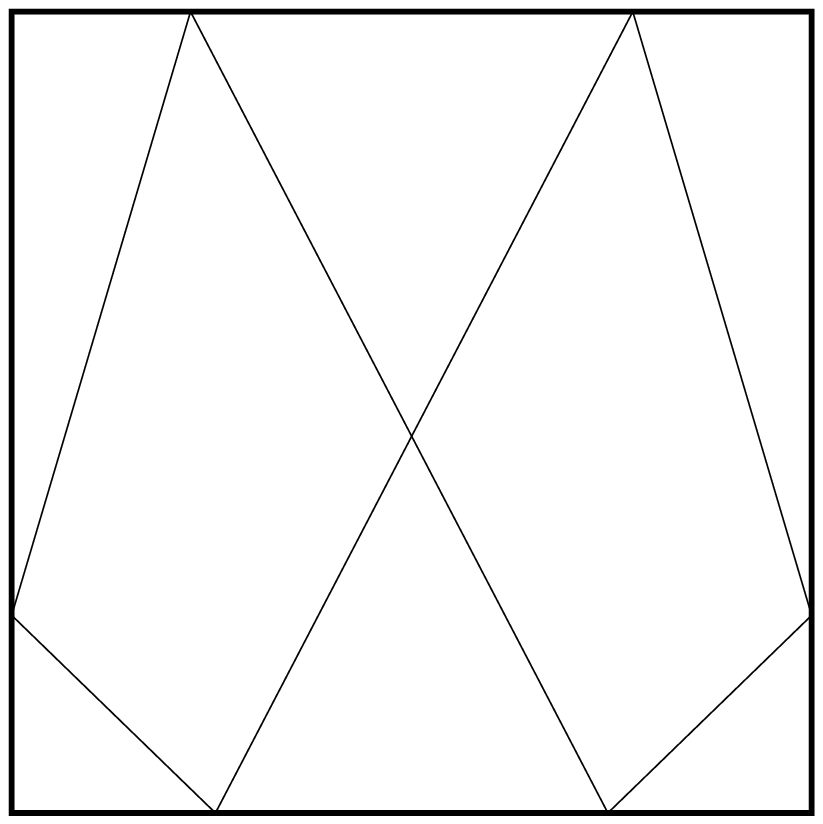}}\quad 
\subfloat[$q_2$]{\includegraphics[width=.2\textwidth]{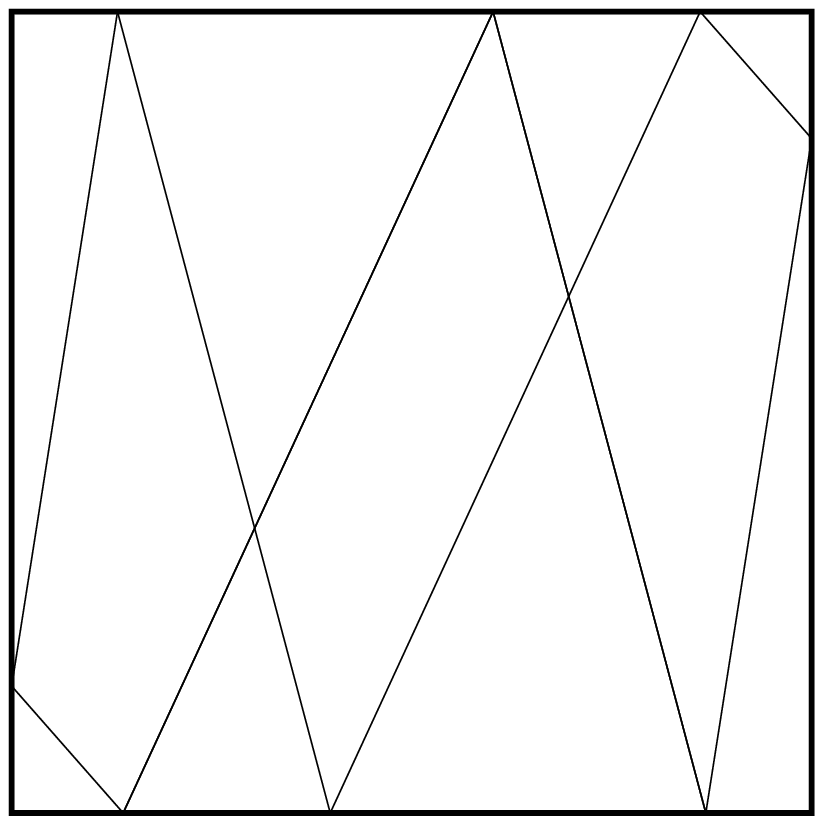}}\quad  
\subfloat[$q_3$]{\includegraphics[width=.2\textwidth]{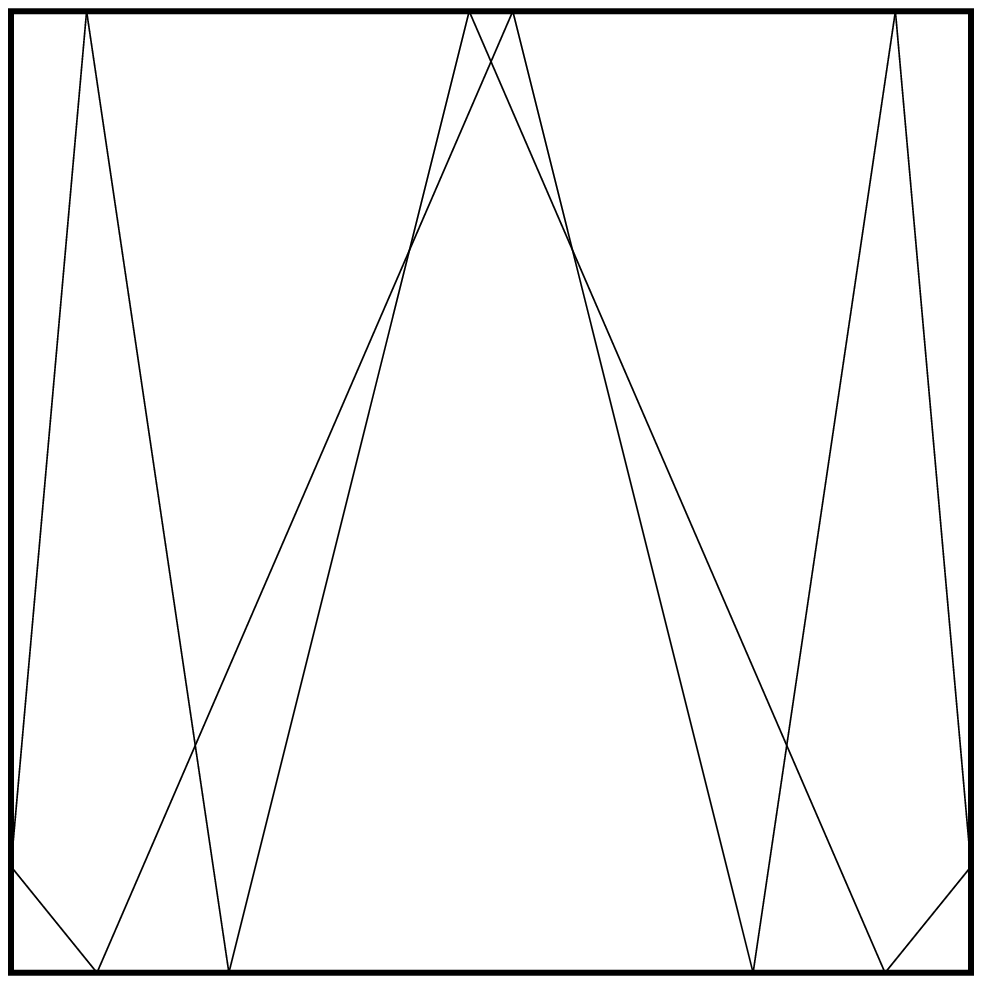}}\\ 
\subfloat[$p_1$]{\includegraphics[width=.2\textwidth]{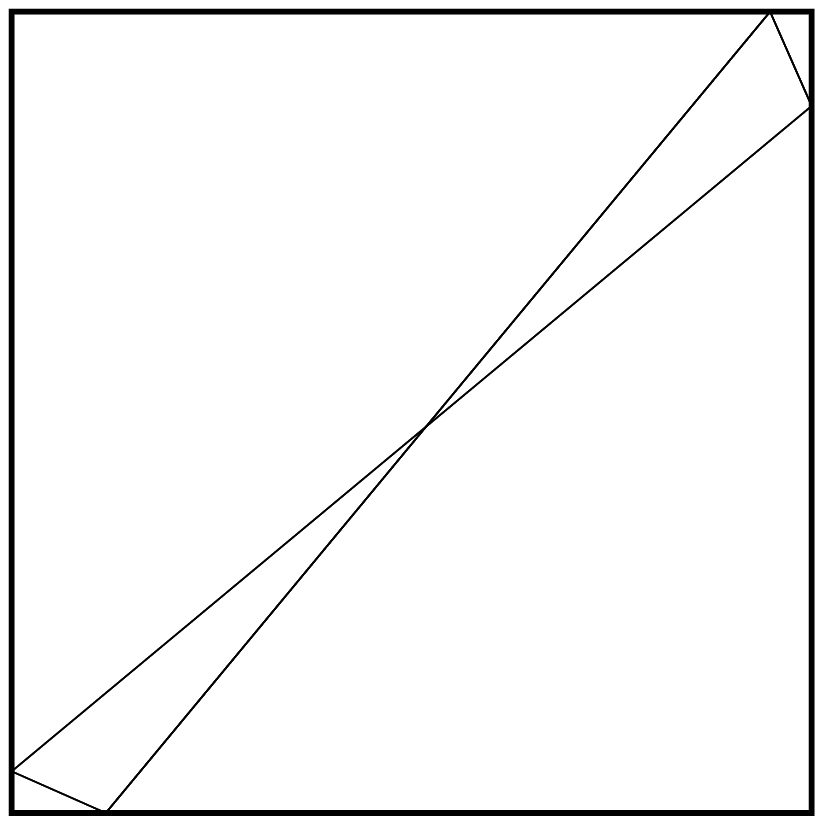}}\quad 
\subfloat[$p_2$]{\includegraphics[width=.2\textwidth]{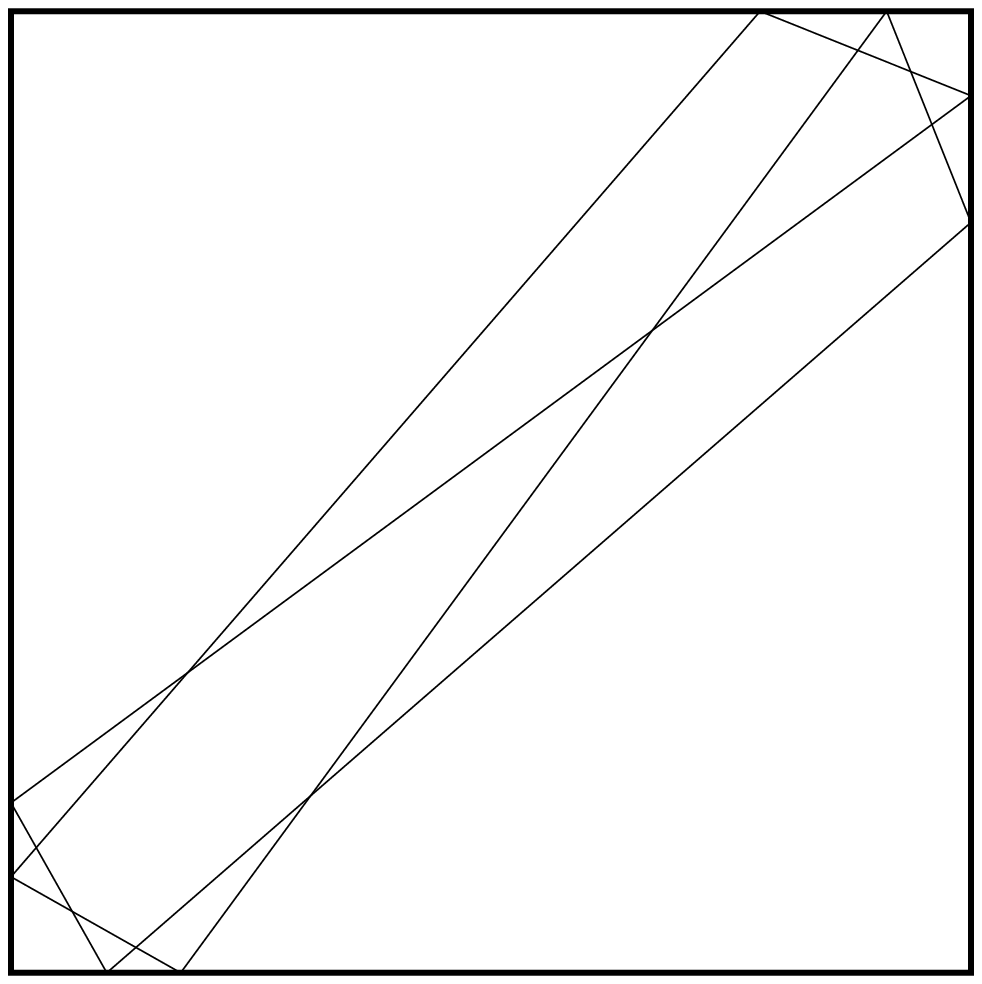}}\quad  
\subfloat[$p_3$]{\includegraphics[width=.2\textwidth]{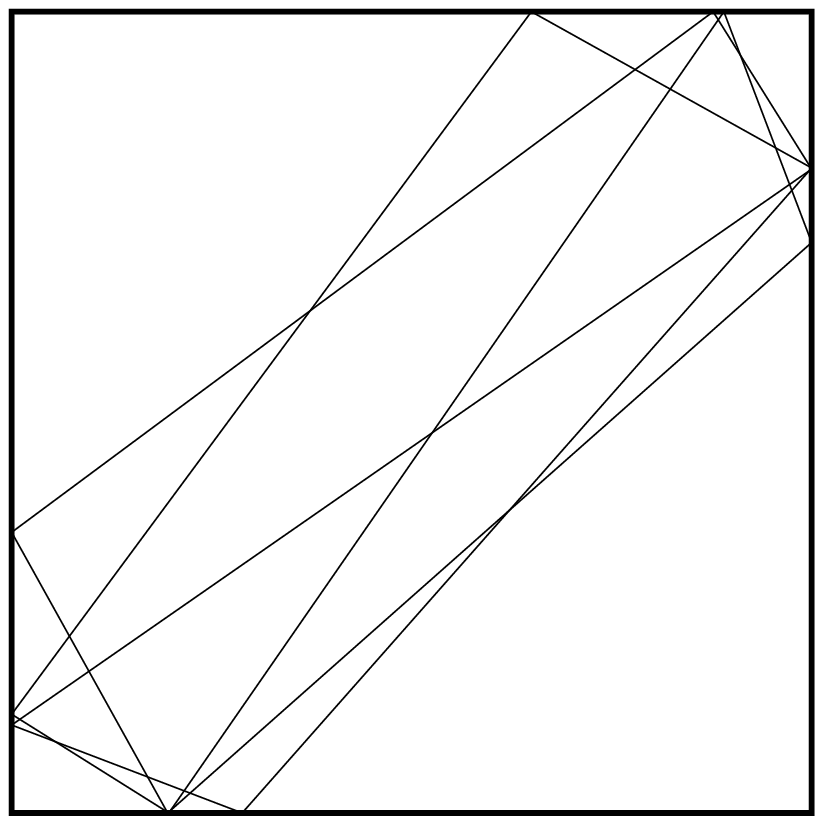}}\\ 
\caption{Billiard trajectories corresponding to the periodic points $q_n$ and $p_n$ for $\lambda=0.6$.} \label{fig periodic orbits}
\end{figure*}

Given $(n,m)\in\mathbb{N}^2$, a straightforward computation shows that
\begin{equation*}
f_2^n\circ f_1^m(s,\theta)=\left((-1)^{n-1}\Upsilon_{n,m}(\theta,s),g_\lambda^n(\lambda^m\theta)\right),
\end{equation*}
where $\Upsilon_{n,m}$ is given by
$$
\Upsilon_{n,m}(\theta,s)=\left[h_{n-1}(\lambda^m\theta)-s-S_{m-1}(\theta)\right]\gamma_{n}(\lambda^m\theta)
$$
and the sequences of functions $h_n$, $\gamma_n$ and $S_n$ are defined by 
\begin{equation*}\begin{split}
h_{n}(\theta)&=\sum_{i=0}^{n}(-1)^i\prod_{j=0}^{i-1}\tan(g_\lambda^j(\theta))\\
\gamma_{n}(\theta)&=\prod_{i=0}^{n-1}\cot(g_\lambda^i(\theta))\\
S_n(\theta)&=\sum_{i=0}^{n}\tan(\lambda^i\theta)\,.
\end{split}
\end{equation*}
Recall that $g_\lambda$ is the affine contraction $g_\lambda(\theta)=\lambda(\pi/2-\theta)$. For each $n\geq 1$, define $q_n$ and $p_n$ as the unique solutions, when they exist, of 
\begin{equation*}
f_2^2\circ f_1^n(q_n)=q_n\quad\text{and}\quad f_2^{2n-1}\circ f_1(p_n)=p_n\,.
\end{equation*}
In agreement with this definition, we set $q_0=p_\lambda$.

\begin{proposition}
	\label{prop:qn}
There exists a unique decreasing sequence $c_n\in(0,\lambda_{2})$ such that $q_n$ is a periodic point of period $n+2$ for $ \phi_\lambda$ if and only if $\lambda\in(0,c_n)$.
\end{proposition}

\begin{proof}
Let $q_n=(s_n,\theta_n)$. A simple computation shows that
	\begin{equation}\label{eq:sn}\begin{split}
	s_n&=(1-S_{n}(\theta_n))\frac{\gamma_{2}(\lambda^n\theta_n)}{\gamma_{2}(\lambda^n\theta_n)-1},\\
	\theta_n&=\frac\pi2\frac{\lambda(1-\lambda)}{1-\lambda^{2+n}}\,.
	\end{split}
	\end{equation}

Let $ \Delta_{n} $ be the set of points $(s,\theta) \in M_1$ such that 
\begin{equation*}
1-S_n(\theta)< s<1-S_{n-1}(\theta)\,.
\end{equation*}
Since $f^{n}_{1}(\Delta_{n}) \subset M_{2}$, to show that $q_{n}$ exists it is enough to check that $(s_{n},\theta_{n})$ belongs to $\Delta_{n}$.

In the following estimates we will frequently use the fact that the tangent is a convex function, i.e., for every $x,y\in(0,\frac\pi2)$ and $0<\lambda<1$ we have $$\tan(\lambda x +(1-\lambda) y)< \lambda\tan(x)+(1-\lambda)\tan(y)\,,$$ and moreover that for every $0<x<\frac\pi2$, $\tan(x)>x$. 

Let us start by proving that $s_n>1-S_{n}(\theta_n)$. By definition
\begin{align*}
\gamma_{2}(\lambda^n\theta_n)&=\cot(\lambda^n\theta_n)\cot\left(\lambda\left(\frac\pi2-\lambda^n\theta_n\right)\right)\\
&>\frac{\cot(\lambda^n\theta_n)}{\lambda\tan\left(\frac\pi2 - \lambda^n\theta_n\right)}\,.
\end{align*}

Thus $\gamma_{2}(\lambda^n\theta_n)>1$ which, taking into account the definition of $s_n$, implies that 
$s_n>1-S_n(\theta_n)$. To prove the other inequality, we start by noting that
\begin{equation*}
s_n=1-S_{n-1}(\theta_n)+\frac{1-S_{n-1}(\theta_n)-\cot(g_\lambda(\lambda^n\theta_n))}{\gamma_2(\lambda^n\theta_n)-1}.
\end{equation*}
Since $\gamma_{2}(\lambda^n\theta_n)>1$, we only need to prove that
\begin{equation}\label{ine:sncot}
S_{n-1}(\theta_n)+\cot(g_\lambda(\lambda^n\theta_n))>1.
\end{equation} 
Using the definition of $S_n$ we get
$$
S_{n-1}(\theta_n)=\sum_{i=0}^{n-1}\tan(\lambda^i\theta_n)>\theta_n\frac{1-\lambda^{n}}{1-\lambda}\,.
$$
On the other hand
$$
\cot(g_\lambda(\lambda^n\theta_n))=\tan(\lambda^{-1}\theta_n)>\frac{\theta_n}{\lambda}\,.
$$
Putting these estimates together we obtain
\begin{align*}
S_{n-1}(\theta_n)+&\cot(g_\lambda(\lambda^n\theta_n))>\theta_n\frac{1-\lambda^{n+1}}{\lambda(1-\lambda)}\\
&=\frac\pi2\frac{1-\lambda^{n+1}}{1-\lambda^{n+2}}>\frac\pi2\frac{n+1}{n+2}>1
\end{align*}
for every $\lambda\in(0,1)$. 

It remains to prove that $s_n>0$. It is clear that $dS_n(\theta_n(\lambda))/d\lambda>0$ for every $\lambda\in(0,1)$. Since $S_n(\theta_n(0))=0$ and $S_n(\theta_n(1))=(n+1)\tan(\pi(n+2)^{-1}/2)>1$ for $n\geq1$, we conclude that $S_n(\theta_n(\lambda))<1$ if and only if $\lambda\in(0,c_n)$. Here $c_n\in(0,1)$ is the unique solution of
$$
\sum_{i=0}^{n}\tan\left(\lambda^i \theta_n(\lambda)\right)=1.
$$
Thus $s_n>1-S_n(\theta_n(\lambda))>0$ if and only if $\lambda\in(0,c_n)$. Now we prove that $\left\{c_{n}\right\}$ is decreasing. Since $dS_n(\theta_n(\lambda))/d\lambda>0$ it is sufficient to prove that $$S_{n+1}(\theta_{n+1})>S_{n}(\theta_n)\,.$$ 
By definition of $S_n$ we have that
\begin{align*}
S_{n+1}(\theta_{n+1})=S_n(\theta_n)&+\sum_{i=0}^{n}\left[\tan(\lambda^i\theta_{n+1})-\tan(\lambda^i\theta_{n})\right]\\
&+\tan(\lambda^{n+1}\theta_{n+1})\,.
\end{align*}
Let $\Delta S_n=S_{n+1}(\theta_{n+1})-S_n(\theta_n)$. Note that $\theta_{n}>\theta_{n+1}$. Since $\tan(y-x)>(\tan(y)-\tan(x))(1-\tan(x)\tan(y))$ for every $0<x<y< \pi/4$ we get
\begin{align*}
\Delta S_n&>\tan(\lambda^{n+1}\theta_{n+1})-\sum_{i=0}^{n}\frac{\tan(\lambda^i(\theta_n-\theta_{n+1}))}{1-\tan(\lambda^i\theta_{n+1})\tan(\lambda^i\theta_{n})}\\
&>\tan(\lambda^{n+1}\theta_{n+1})-\sum_{i=0}^{n}\frac{\tan(\lambda^i(\theta_n-\theta_{n+1}))}{1-\tan^2(\theta_{n})}\,.
\end{align*}
In the derivation of the previous inequality we have used the upperbound: $\tan(\lambda^i\theta_{n+1})\tan(\lambda^i\theta_{n})<\tan^2(\theta_{n})$ for every $n\geq1$.
Let 
$$
\rho_n=\frac{\lambda^{n+2}}{1+\lambda+\cdots+\lambda^{n+1}}.
$$
Clearly $\theta_n=(1+\rho_n)\theta_{n+1}$. Since $\rho_n/\lambda^{n+1}<1$ we obtain
\begin{align*}
\Delta S_n&>\tan(\lambda^{n+1}\theta_{n+1})-\sum_{i=0}^{n}\frac{\tan(\rho_n\lambda^i\theta_{n+1})}{1-\tan^2(\theta_{n})}\\
&>\tan(\lambda^{n+1}\theta_{n+1})\left(1-\frac{\rho_n(1-\lambda^{n+1})}{\lambda^{n+1}(1-\tan^2(\theta_{n}))(1-\lambda)}\right)\,.
\end{align*}
Using the expression for $\rho_n$ we get
$$
\Delta S_n>\frac{\tan(\lambda^{n+1}\theta_{n+1})(1-\lambda)}{(1-\tan^2(\theta_{n}))(1-\lambda^{n+2})}\left(1-\frac{\pi \lambda\tan^2(\theta_n)}{2\theta_n}\right)\,.
$$
Since $0<\theta_n<\theta_1<\pi/4$ we have that
\begin{align*}
\tan^2(\theta_n)&<\left(\frac{\theta_n}{\theta_1}\right)^2\tan^2(\theta_1)\\
&<\frac{2 \theta_n}{\pi\lambda}(1+\lambda+\lambda^2)\tan^2\left(\frac{\pi}{2\lambda(1+\lambda+\lambda^2)}\right)\\
&<\frac{2 \theta_n}{\pi \lambda}3\tan^2\left(\frac\pi6\right)=\frac{2 \theta_n}{\pi\lambda}\,.
\end{align*}
Hence $\Delta S_n>0$ as we wanted to show. Finally, $\left\{c_n\right\}$ is bounded from above by $\lambda_2$ since a numerical computation reveals that
$$
c_1=0.7964\ldots
$$
which, taking into account the numerical value of $\lambda_2$, implies that $c_n<c_1<\lambda_2$.
\end{proof}

The proof of the next result is omitted because it is similar to the previous one.
\begin{proposition}
If $\lambda\in(0,\lambda_2)$ then $p_n$ is a periodic point of period $2n$ for $\phi_\lambda$.
\end{proposition}

By Lemma~\ref{le:periodic-points}, these periodic points are hyperbolic. As we shall see in the next section, these orbits seem to play an important role in the dynamics of $\phi_\lambda$ for different values of $\lambda$. The corresponding billiard orbits in configuration space are depicted in Fig.~\ref{fig periodic orbits}.

Since the sequence $ c_{n} $ is decreasing, we can define $\lambda_1=\lim_{n\to\infty} c_n$. The number $\lambda_{1}$ is also the unique solution of the equation
$$
\sum_{n=0}^{\infty}\tan\left(\frac\pi2\lambda^{i+1}(1-\lambda)\right)=1 \qquad \text{for } \lambda \in (0,1).
$$
In geometrical terms, when $\lambda=\lambda_1$, the intersection point of the curve $S_\infty$ with the line $\theta=\pi\lambda(1-\lambda)/2$ lies exactly on the line $s=0$ (see Fig.~\ref{trapping:region}(a)). This intersection point is also the limit of the sequence $q_n$. A numerical computation shows that
\[
\lambda_{1}=0.6218...
\]

By Proposition \ref{prop:qn}, the periodic points $q_n$ disappear as $\lambda$ increases from $\lambda_1$ to $\lambda_2$. The point $q_1$ is the last to disappear for a value of $\lambda$ close to $\lambda_2$. All points $ q_{n}$ are contained in the light-colored trapping region depicted in Fig.~\ref{trapping:region}(a).

\begin{proposition} \label{homcl:qn}
We have $W^u(q_n)\cap B\neq \emptyset$ for every $0<\lambda<\lambda_1$ and every $n$ sufficiently large. In particular, the homoclinic class of $q_n$ is a transitive horseshoe provided that $n$ is sufficiently large.
\end{proposition}
\begin{proof}
Let $q_n=(s_n,\theta_n)$ where $s_n$ and $\theta_n$ are given by \eqref{eq:sn}. Hence as $n\to\infty$  these points approach the horizontal line $\theta=\pi\lambda(1-\lambda)/2$. It is also straightforward to check that $W^ u(q_n)$ contains the horizontal segment joining $(0,\theta_n)$ to $q_n$. Let 
$$
\sigma(\theta)=1-\sum_{i=0}^\infty \tan(\lambda^i \theta)
$$
be the map whose graph is $S_\infty$, the upper bound of $B$. This function is decreasing in $\lambda$, and by definition of $\lambda_1$ we have $\sigma(\pi\lambda(1-\lambda)/2)=0$ when $\lambda=\lambda_1$. Also, we must have $\sigma(\pi\lambda(1-\lambda)/2)>0$ for every $0<\lambda<\lambda_1$. Thus, $\sigma(\theta_n)>0$ for all large enough $n$. This proves that $W^ u(q_n)$, which contains the segment joining $(0,\theta_n)$ to $q_n$, intersects the set $B\subset W^ s(P)$, bounded from above by the graph of $\sigma$. Hence this homoclinic class is a horseshoe. To complete the proof, we just need to observe that transitivity is a general property of the homoclinic classes (see Ref.~\cite[Ch.~IX, Prop.~5.2]{r}).
\end{proof}

\begin{proposition}\label{prop:Gamma}
The following statements hold for $\lambda>\lambda_1$:
\begin{enumerate}
	\item $W^s(P)=B\cup\phi^{-1}_\lambda(B)$.
	\item $\Gamma=M^{+}\setminus W^s(P)$ is a trapping region, i.e. $\phi_\lambda(\Gamma)\subset \mathrm{int}(\Gamma)$.
	\item $W^u(q_n)\cap B= \emptyset$ for every $n\geq1$.
\end{enumerate}
\end{proposition}
\begin{proof}
Since the second and third statements immediately follow from the first one, we only prove that $$
W^s(P)=B\cup\phi^{-1}_\lambda(B)\,.
$$ By definition, $W^s(P)=\bigcup_{n\geq0}\phi_\lambda^{-n}(B)$. Thus, 
$$
B\cup\phi^{-1}_\lambda(B)\subseteq W^s(P)\,.
$$ 
To prove the opposite inclusion it is sufficient to prove that
$$
\phi^{-2}_\lambda(B)\subset B\cup\phi^{-1}_\lambda(B)\,.
$$
Suppose that the previous inclusion does not holds, i.e. there exists $x\in\phi^{-2}_\lambda(B)$ such that neither $x\in B$ nor $x\in\phi^{-1}_\lambda(B)$. Thus $\phi_\lambda(x)\in\phi^{-1}_\lambda(B)\setminus B$. On the other hand, it is clear that for every $\lambda>\lambda_1$ we have
$$
(\phi^{-1}_\lambda(B)\setminus B) \cap \phi_\lambda(M)=\emptyset\,,
$$
yielding a contradiction. 
\end{proof}

\begin{figure*}
\subfloat[$\lambda=0.615\in(\lambda_0,\lambda_1)$]{\includegraphics[width=.3\textwidth]{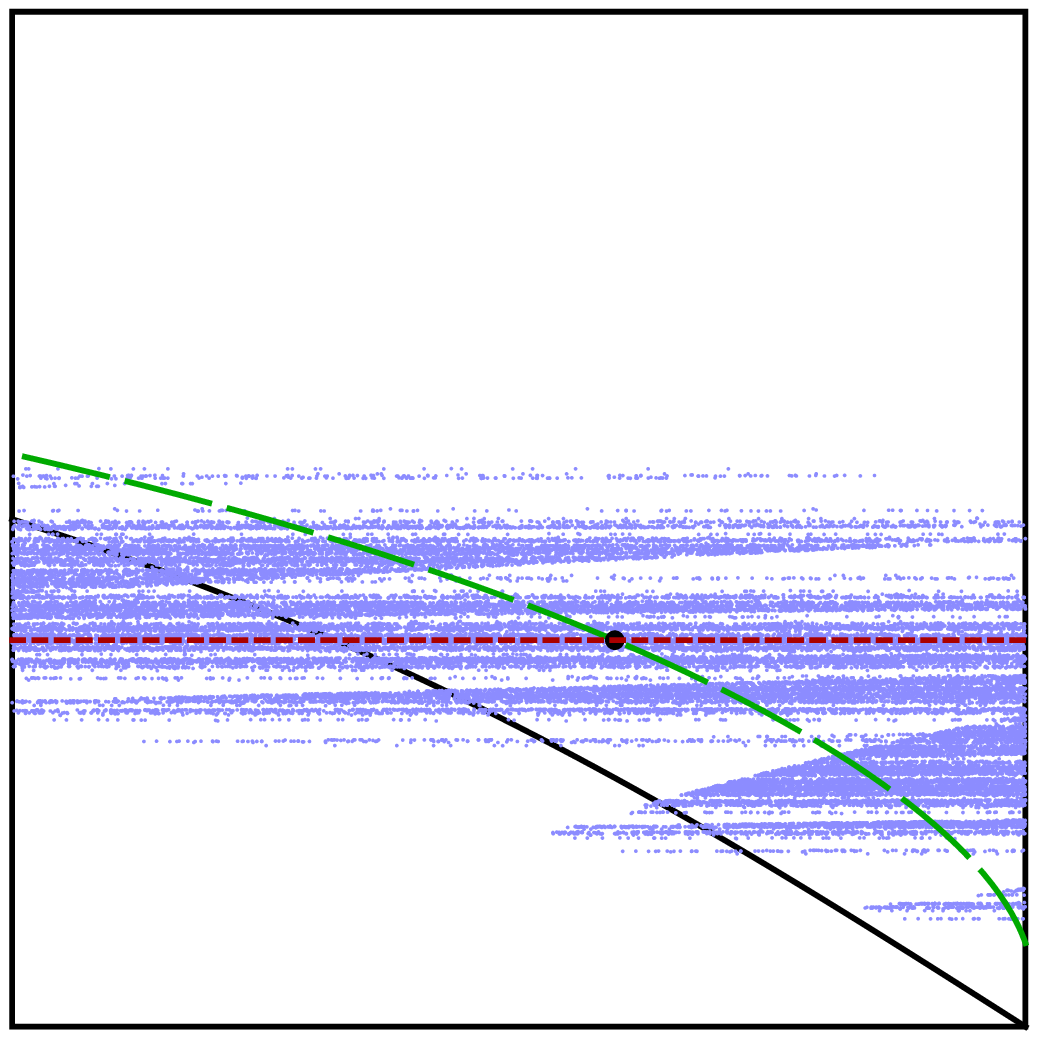}} 
\subfloat[$\lambda=0.75\in(\lambda_1,\lambda_2)$]{\includegraphics[width=.3\textwidth]{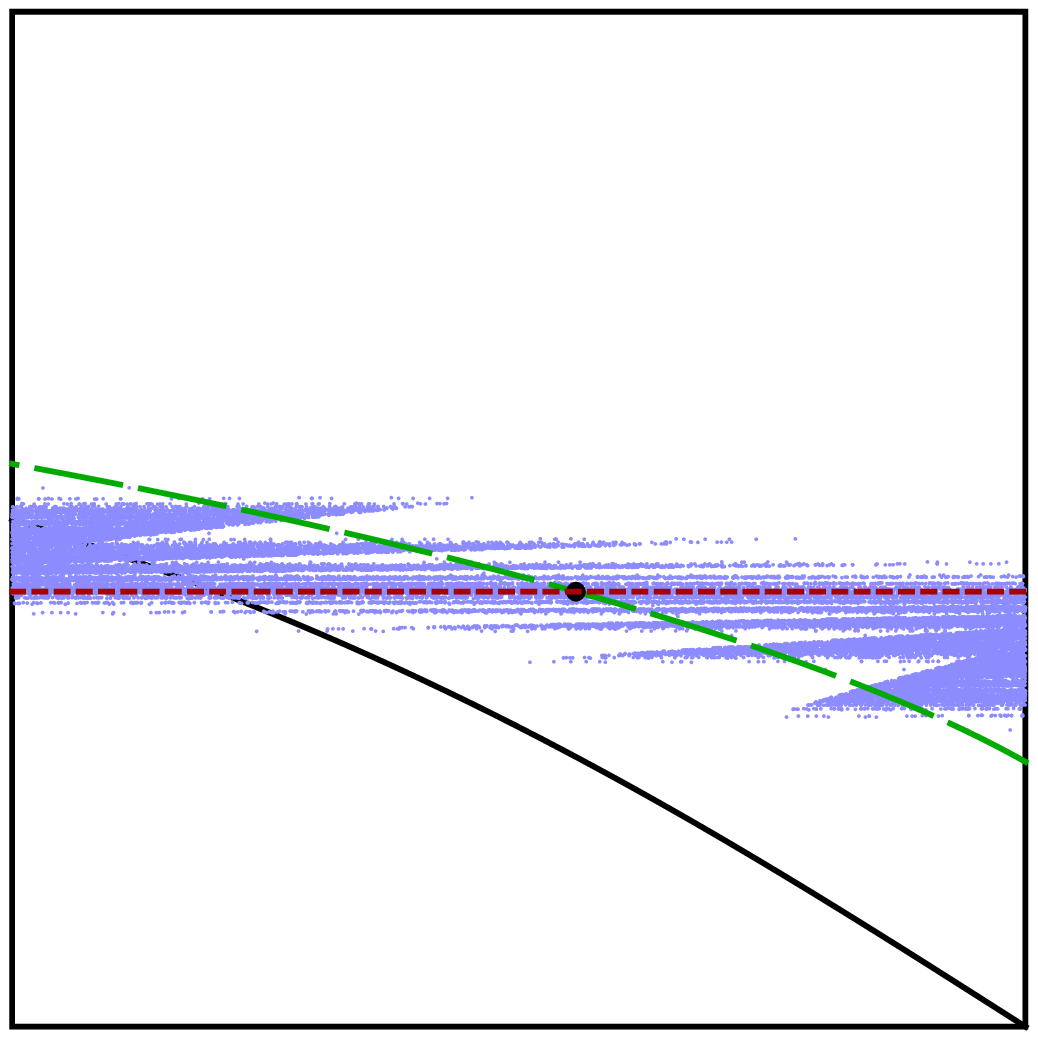}}  
\subfloat[$\lambda=0.88\in(\lambda_2,1)$]{\includegraphics[width=.3\textwidth]{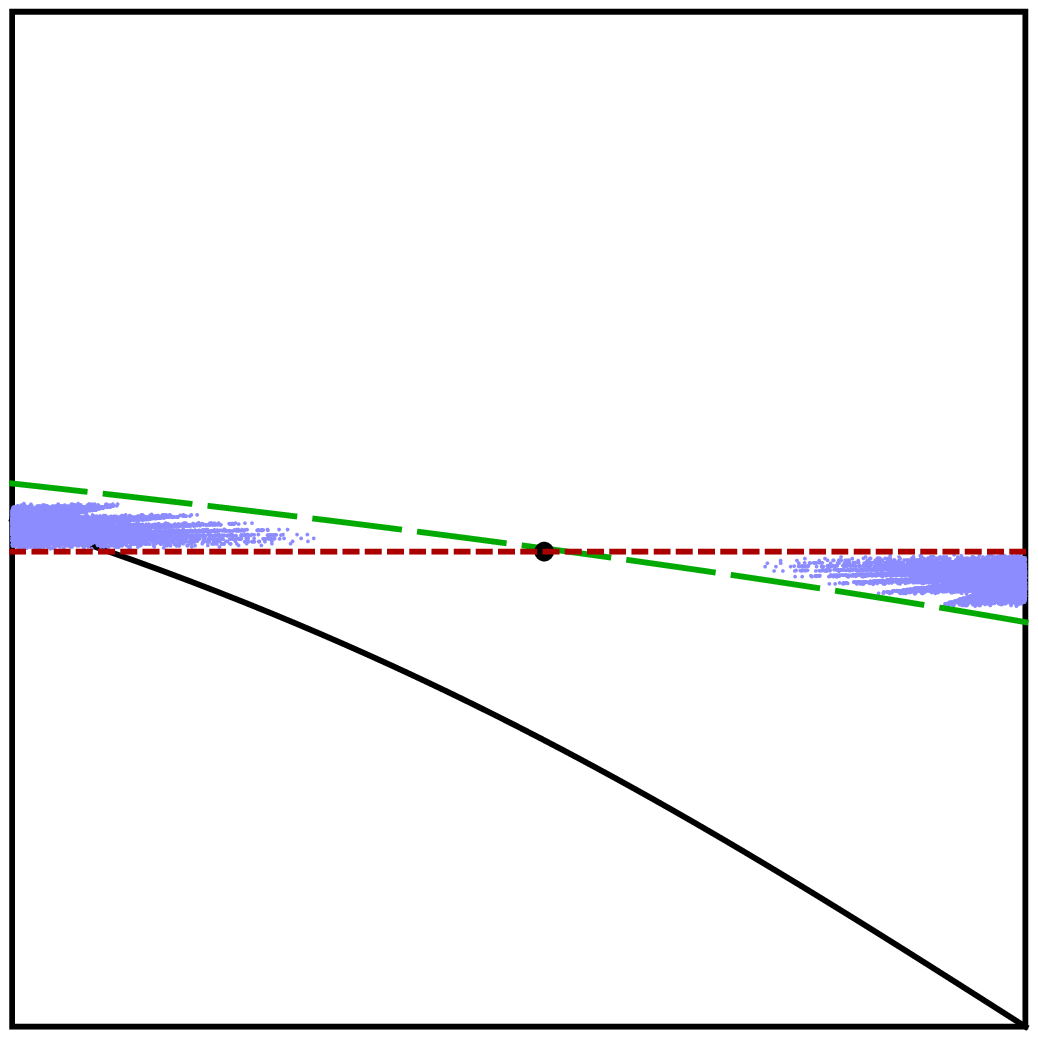}} 
\caption{Local stable (dashed green curve) and unstable (dotted red curve) manifolds of $ p_{\lambda} $, and attractor $A_\lambda$ (blue region).}
\label{figures_attractor}
\end{figure*}

\begin{proposition}\label{p:qn}
There exists $\delta>0$ such that the periodic points $q_n$ are all homoclinically related with $p_\lambda$ for every $0<\lambda<\delta$.
\end{proposition}
\begin{proof}
For $\lambda>0$ close to $0$, the local  stable manifold $W^s_{\mathrm{loc}}(p_\lambda)$ is the graph of the concave monotonic function $h_\lambda$ connecting the left side $s=0$ to the vertex $(s,\theta)=(1,0)$. Using a graph transform argument we can prove that the
local  stable manifold $W^s_{\mathrm{loc}}(q_n)$ is the graph of a
concave monotonic function $s=h_n(\theta)$ with $\sigma_{n+1}(\theta)< h_n(\theta) < \sigma_{n}(\theta)$,
where 
$$\sigma_n(\theta)  =  1-\left( \tan(\theta)+\tan(\lambda \theta) + \ldots +\tan(\lambda^{n-1}\theta)\right)\;.$$
Notice that $M_{1}^ n=\{\,(s,\theta)\, :\, \sigma_{n+1}(\theta)< s < \sigma_{n}(\theta)\,\}$ is the region of all points in $M_1$ mapped by $f_1^ n$ into the domain $M_2$. The graph $s=h_n(\theta)$ also connects the left side $s=0$ to the vertex $(s,\theta)=(1,0)$. An easy computation shows that 
$$
1-\sum_{i=1}^{\infty}\tan(\lambda^i\theta_n)>0
$$ 
and $\theta_n> \pi \lambda/4$, the second inequality for $\lambda<1/2$. Because $f_1(0,\pi/2)=(1,\pi \lambda/4)$ we can deduce from the inequality $\theta_n>\pi \lambda/4$ that the local unstable manifold of $q_n$ is the horizontal segment connecting $(0,\theta_n)$ to $(1,\theta_n)$. Whence $W^u_{\rm loc}(q_n)=[0,1]\times \{\theta_n\}$ intersects $W^s_{\rm loc}(q_m)={\rm graph}(h_m)$ for any pair of integers $n,m\geq 1$. Since $p_\lambda=q_0$, the proposition is proved.
\end{proof}

For future use, we now introduce the new constant
$$ \lambda_0=\inf\{ \lambda>0 \,:\, \exists\, n\geq 1\,,\, W^u(p_\lambda)\cap W^ s(q_n)=\emptyset\,\} \;. $$
From Proposition~\ref{p:qn}, it follows that $\lambda_0 \ge \delta > 0$, where $ \delta $ is as in Proposition~\ref{p:qn}. Numerically, we found that
$$\lambda_{0} > 0.6104\;.$$ 


\section{Bifurcation of the limit set}
\label{se:bifurcation}

Recall that $\Omega_{\lambda}$ is the nonwandering set of the map $\phi_\lambda$. In this last section, we formulate a conjecture on the decomposition of $\Omega_{\lambda}$, and discuss the changes in this decomposition as the parameter $\lambda$ varies.

\begin{conjecture}\label{main:conjecture}
For any $0<\lambda<1$, the nonwandering set $ \Omega_{\lambda} $ is a union of three sets:
$$
\Omega_{\lambda}=P\cup H_\lambda\cup A_\lambda,
$$
where $P$ is the parabolic attractor introduced in Section \ref{ss:parabolic_attractor}, $A_\lambda$ is a hyperbolic transitive attractor, and $H_\lambda$ is a horseshoe. Moreover, $ H_{\lambda} $ is either transitive or else a (possibly empty) hyperbolic chain of transitive horseshoes. In particular,  
$$
M^{+}=W^s(P)\cup W^s(H_\lambda)\cup W^s(A_\lambda).
$$
\end{conjecture}

Our next conjecture is justified by the fact that $P$ consists of periodic points, and the set $A_{\lambda} \cup H_{\lambda}$ is hyperbolic. 

\begin{conjecture}\label{conj:per} The set of periodic points is dense in $\Omega_{\lambda}$.
\end{conjecture}


The rest of the section is devoted to the justification of the previous Conjecture~\ref{main:conjecture}, and to the analysis of the changes in the sets $ H_{\lambda} $ and $ A_{\lambda} $ as $ \lambda $ varies. The conclusions based on numerical observations are presented as conjectures, whereas the conclusions based on analytical arguments are presented as propositions with their proofs. We split our discussion into four parts, each corresponding to $\lambda$ taking values inside one of the following intervals: $ (0,\lambda_{0}) $, $ (\lambda_{0},\lambda_{1}) $, $ (\lambda_{1},\lambda_{2}) $ and $ (\lambda_{2},1)$. See Fig.~\ref{figures_attractor}.

\begin{figure*}
\subfloat[]{\includegraphics[width=0.3\textwidth]{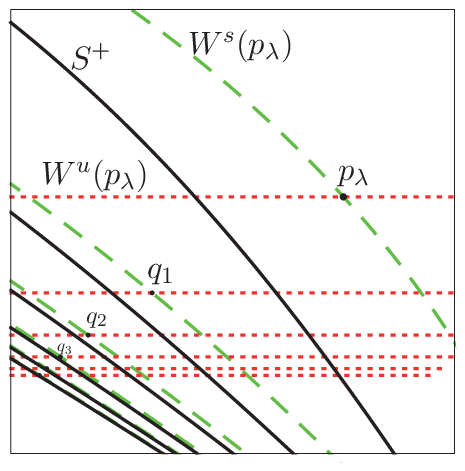}}             
\hspace{3.5cm}
\subfloat[]{\includegraphics[width=0.3\textwidth]{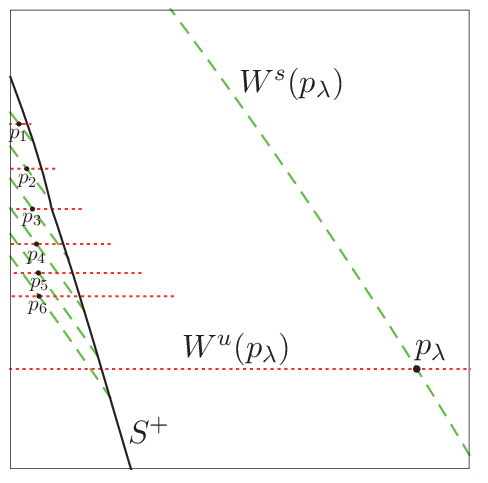}}
\caption{Zoom of the phase space with maximal local invariant manifolds of the periodic points. (a) Points $q_n$ together with their local stable (dashed green) and unstable (dotted red) manifolds for $\lambda=0.6$. The black curves represent some iterates of the singular set $S^{+}$. \\ \hspace{1.35 cm} (b) Points $p_n$ together with their  local stable (dashed green) and unstable (dotted red) manifolds for $\lambda=0.85$.}
\label{figure:hom:related}
\end{figure*}

\subsection {  $0<\lambda<\lambda_0$ }

The following conjecture is suggested by numerical computations of the invariant manifolds for the points $q_n$ (see Fig.~\ref{figure:hom:related}(a)).

\begin{conjecture}\label{qn:qn+1}
The manifolds $W^u(q_{n+1})$ and $W^s(q_n)$ intersect transversally for every $0<\lambda<c_{n+1}$. Moreover, provided that $ 0<\lambda<\lambda_{1} $, all $q_n$ are mutually homoclinically related for sufficiently large $ n $.
\end{conjecture}

In light of Conjecture~\ref{conj:per}, the next conjecture simply states that
$\Omega_{\lambda}\setminus P$ is the union of the homoclinic classes of the $q_{n}$.

\begin{conjecture}\label{qn:x}
Suppose that $ x \in \Omega_{\lambda} \setminus P $ is a periodic point of $ \phi_{\lambda} $. Then there are $ m,n \ge 0 $ such that $W^u(x)$ and $W^s(q_n)$ intersect transversally, and $W^u(q_m)$ and $W^s(x)$ intersect transversally.
\end{conjecture}

\begin{proposition} \label{horse:shoe}
If Conjectures~\ref{qn:qn+1} and \ref{qn:x} hold, then $\Omega_{\lambda}=P\cup H_\lambda$ and $H_\lambda$ is a transitive horseshoe for $0<\lambda<\lambda_0$.
\end{proposition}
\begin{proof}
By definition of $\lambda_0$ all $q_n$ are homoclinically related for  $0<\lambda<\lambda_0$. By Conjectures~\ref{qn:qn+1} and \ref{qn:x}, $\Omega_{\lambda}=P\cup H_\lambda$, where $H_\lambda$ denotes the homoclinic class of $p_\lambda=q_0$. Then Proposition~\ref{homcl:qn} shows that $H_\lambda$ is a transitive horseshoe.
\end{proof}

\subsection {$\lambda_0<\lambda<\lambda_1$}

In this parameter range the set $\Omega_{\lambda}\setminus P$ splits into two or more homoclinic classes dynamically partially ordered. At the bottom of this hierarchy of homoclinic classes lies a transitive hyperbolic attractor, and at the top a transitive horseshoe whose unstable set intersects the basin of attraction of $P$. We write $H\prec H'$ for $W^u(H')\cap W^ s(H)\neq \emptyset$.

\begin{proposition} \label{Cis}
If Conjectures~\ref{conj:per}-\ref{qn:x} hold, then for $\lambda_0<\lambda<\lambda_1$, there exists $N \ge 1 $ such that $\Omega_{\lambda}=P\cup C_0\cup C_1\cup \ldots \cup C_N$ and
\begin{enumerate}
\item each $C_i$ is the homoclinic class of some periodic point $q_{n_i}$,
\item $C_i\cap C_j=\emptyset$ whenever  $i\neq j$,
\item $C_0\prec C_1\prec \ldots \prec C_N$,
\item $C_0$ is a transitive hyperbolic attractor,
\item $P\prec C_N$.
\end{enumerate}
\end{proposition}
\begin{proof}
Let $C_0, C_1,\ldots, C_N$ be the homoclinic classes of the periodic points $q_n$.
By conjecture~ \ref{qn:x}, we have $\Omega_{\lambda} \setminus P=C_0\cup C_1\cup\ldots \cup C_N$.
These sets are obviously disjoint.
Conjecture~\ref{qn:qn+1} implies the sets $C_i$ are ordered in a finite chain,
 and we can always display them as in item 3.
$C_0$ is the homoclinic class  of the fixed point $q_0=p_\lambda$, and hence
a transitive invariant set. It is  attracting since it lies at the chain's bottom, and it is hyperbolic because of Proposition \ref{pr:hyperbolicity}. The set $H=C_1\cup \ldots \cup C_N$ is a chain of transitive hyperbolic horseshoes. Finally, since $C_N$ is at the chain's top, Proposition~\ref{homcl:qn} implies that $P\prec C_N$.
\end{proof}

By the definition of $\lambda_0$, for every $\lambda_0<\lambda<\lambda_1$ there is some $n\geq 1$  such that $W^u(q_0)\cap W^s(q_{n+1})=\emptyset$, and, in view of Conjecture~ \ref{qn:qn+1}, this implies there is some $n\geq 1$ such that $W^u(q_{n})\cap W^s(q_{n+1})=\emptyset$. Given $n\geq 1$, let $\overline{\lambda}_n$ be the bifurcation point where the homoclinic connection $W^u(q_{n})\cap W^s(q_{n+1})$ breaks down. The numerical value given above for $\lambda_0$ was obtained from the following dichotomy: for $\lambda<\lambda_0$ almost every point is attracted to $P$, while for $\lambda >\lambda_0$ there is a non trivial hyperbolic attractor with an open
basin of attraction. We did not try to understand these heteroclinic connection breaking bifurcations $\overline{\lambda}_n$, but numerical plots indicate that $\lambda_0=\overline{\lambda}_n$ for some rather small $n$, probably $n\leq 3$.

\subsection {  $\lambda_1<\lambda<\lambda_2$ }

In this parameter range the periodic points $q_n$ vanish one by one.
More precisely, according to Proposition~\ref{prop:qn} there is a decreasing sequence of bifurcation parameters
$$\lambda_1 <\ldots <c_{n+1} < c_n <\ldots < c_2 < c_1 < \lambda_2,$$
and $q_n$ persists for $\lambda <c_n$, but vanishes for $\lambda>c_n$.
Hence, unlike the previous interval, only finitely many $q_n$ persist for each $\lambda_1<\lambda<\lambda_2$.

\begin{proposition} \label{Cns}
If Conjectures~\ref{conj:per}-\ref{qn:x} hold, then for $\lambda_1<\lambda<\lambda_2$, there exists $N \ge 0$ such that $\Omega_{\lambda}=P\cup C_0\cup C_1\cup \ldots \cup C_N$ and 
\begin{enumerate}
\item each $C_i$ is the homoclinic class of some periodic point $q_{n_i}$,
\item $C_i\cap C_j=\emptyset$ whenever  $i\neq j$,
\item $C_0\prec C_1\prec \ldots \prec C_N$,
\item $C_0$ is a transitive hyperbolic attractor,
\item $P \nprec C_i$ for all $i=0,1,\ldots, N$.
\end{enumerate}
\end{proposition}
\begin{proof}
Keeping the notation of last section, the proof here is a simple adaptation of that of Proposition~ \ref{Cis}.
As before,  the hyperbolic attractor $C_0$ is  the homoclinic class of the fixed point $q_0=p_\lambda$. The main difference is that for $\lambda>\lambda_1$, by Proposition~\ref{prop:Gamma} we have $W^ s(P)=B\cup \phi_\lambda^{-1}(B)$ and there is a trapping region $\Gamma$, disjoint from $W^ s(P)$, forward invariant under $\phi_\lambda$,
which contains all periodic points $q_n$. This proves item 5.
\end{proof}

We found numerically that
\begin{itemize}
\item  $\Omega_{\lambda}\setminus P=C_0$, for $c_1<\lambda<\lambda_2$, 
\item  $\Omega_{\lambda}\setminus P=C_0\cup \{q_1\}$\, with \, $q_1 \nprec q_0$, for $c_2<\lambda<c_1$, 
\item $\Omega_{\lambda}\setminus P=C_0\cup \{q_1\}\cup\{q_2\}$  with   $q_1 \nprec q_0$  
and   $q_2 \nprec q_1$, for $c_3<\lambda<c_2$.
\end{itemize}

\subsection { $\lambda_2<\lambda<1$}
By Proposition~\ref{prop:delta}, in this parameter range the shadowed region $\Delta$ in Fig.~\ref{trapping:region}(b) is a trapping region, i.e. $\Delta$ is forward invariant under $\phi_\lambda$. Moreover, all periodic points $p_n$ must lie inside $\Delta$ whenever they exist.

Our numerical analysis suggests the following conjectures (see Fig.~\ref{figure:hom:related}(b)).  

\begin{conjecture}\label{pn:17}
The periodic points $p_n$ with $n\leq 16$ persist for $\lambda_2 < \lambda < 1$,
while those with $ n \ge 17 $ persist for $\lambda\in (0,a_n]\cup [b_n,1)$,
where the sequences bounding the gap satisfy  $a_n \searrow \lambda_2$ and $b_n \nearrow 1$. In particular, for any given $\lambda_2<\lambda<1$, only finitely many points $p_n$ persist.
\end{conjecture}
 
\begin{conjecture}\label{pn:two} The periodic points $p_n$ generate two homoclinic classes
\begin{enumerate}
\item $C_0$ the homoclinic class of the  $p_n$ with $n\leq 16$ or $\lambda <a_n$,
\item $C_1$ the homoclinic class of the  $p_n$ with $n\geq 17$ and $\lambda >b_n$.
\end{enumerate}
For $\lambda>b_{17}$ (i.e. when $C_1$ becomes non-empty), $C_0\prec C_1$.
\end{conjecture}

The next conjecture simply states that $\Omega_{\lambda}\setminus P$ is the union of the homoclinic classes of the  periodic points $p_n$ and the fixed point $p_\lambda$.

\begin{conjecture}\label{pn:x}
For every  $\lambda_2<\lambda<1$ and every periodic point $x\in \Omega_{\lambda}\setminus (P\cap\{p_\lambda\})$, there exist $n,m\geq 1$ such that $W^u(x)$ and $W^s(p_n)$ intersect transversally, and $W^u(p_m)$ and $W^s(x)$ intersect transversally .
\end{conjecture}

The proof of the following proposition is similar to Proposition~\ref{Cis}.

\begin{proposition}
If Conjecture~\ref{conj:per} and Conjectures~\ref{pn:17}-\ref{pn:x} hold, then for every  $\lambda_2<\lambda<1$, we have $\Omega_{\lambda}=P\cup C_0\cup C_1 \cup \{p_\lambda\}$ and
\begin{enumerate}
\item $C_0$ is a transitive hyperbolic attractor,
\item $C_1$ is a transitive horseshoe (possibly empty),
\item $C_0\cap C_1=\emptyset$,
\item $C_0 \prec C_1\prec \{p_\lambda\}$.
\end{enumerate}
\end{proposition}
 
\begin{acknowledgments}
The authors were supported by Funda\c c\~ao para a Ci\^encia e a Tecnologia through the Program POCI 2010 and the Project ``Randomness in Deterministic Dynamical Systems and Applications'' (PTDC-MAT-105448-2008). G. Del Magno would like to thank M.~Lenci and R.~Markarian for useful discussions.
\end{acknowledgments}



\end{document}